# A Special Issue on Statistical Challenges and Opportunities in Electronic Commerce Research

**Wolfgang Jank and Galit Shmueli**

This special issue is a product of the *First Interdisciplinary Symposium on Statistical Challenges and Opportunities in Electronic Commerce Research*, which took place on May 22–23, 2005, at the Robert H. Smith School of Business, University of Maryland, College Park (www.smith.umd.edu/dit/statschallenges/). The symposium brought together, for the first time, researchers from statistics, information systems, and related fields, all of whom work or are interested in empirical research related to electronic commerce. The goal of the symposium was to cross the borders, discuss joint research opportunities, expose this field and its statistical challenges, and promote collaboration between the different fields.

## 1. BACKGROUND

Electronic commerce (eCommerce) has experienced an extreme surge of popularity in recent years. It has had a huge impact on the way we live today compared to a decade or so ago: It transformed the economy, eliminated borders, opened the door to many innovations and created new ways in which consumers and businesses interact. Although many predicted the death of eCommerce with the "burst of the Internet bubble" in the late 1990s, today eCommerce is thriving more than ever.

By eCommerce we mean any forms of electronic transactions that are commerce-related: online buying, selling or investing; electronic marketplaces like www.amazon.com or online auctions like www.ebay.com; clickstream data and cookie-tracking; e-bookstores and e-grocers; web-based reservation systems and ticket purchasing; marketing email and message postings on web-logs; downloads of music, video and other information; user groups and electronic communities; online discussion boards and learning facilities; open source projects; online banking, and many, many more.

The public nature of many Internet transactions created new opportunities for researchers to gather and analyze data in order to learn about individuals, companies and societies. Theoretical results, founded in economic and game-theoretic models, and derived for the offline, brick-and-mortar world, have often proven not to hold in the online environment. Possible reasons are the worldwide reach of the Internet, its unlimited resources, constant availability and continuous change. For this reason, and also due to the availability of massive amounts of publicly available high-quality web-data, empirical research is thriving.

To date, the fast growing empirical eCommerce research has been concentrated almost entirely in the fields of information systems, economics and marketing. However, as we found out in collaborative work with colleagues in these fields, eCommerce data arrive with many new statistical challenges, ranging from data collection to data exploration, and to analysis and modeling—all of which have generally been overlooked in the current literature. The most likely reason for this is the absence of statisticians from this field. The question is then "*Where are the statisticians?*" We believe that the explanation for the current absence lies in (1) the physical disconnect between academic information systems departments (that are usually housed in business schools) and statistics departments (that tend to be in schools of social sciences, engineering or liberal arts and sciences), and (2) the format in which eCommerce data typically arrive. The first question we hear from statisticians who see our work is "Where did you get the data?" eCommerce data are often organized in HTML pages on the web, which is different from the way our discipline traditionally receives data. For instance, bids placed during an online auction are recorded as a sequence of numbers displayed on an eBay HTML page. Similarly, postings on web-logs ("blogs") are displayed as text on HTML pages. Collecting HTML pages and migrating them into standard databases and spreadsheets is usually done electronically using software programs called "web-agents." However, writing web-agents is not a skill that most statisticians are familiar with, and is not part of the current statistics curriculum. Thus, technology (or the lack thereof) has blinded the eye of the statistics community to an important, huge, rich and high-quality data source.

*Wolfgang Jank is Assistant Professor of Management Science and Statistics, Department of Decision and Information Technologies, Robert H. Smith School of Business, University of Maryland, College Park, Maryland 20742, USA e-mail: wjank@rhsmith.umd.edu. Galit Shmueli is Assistant Professor of Management Science and Statistics, Department of Decision and Information Technologies, Robert H. Smith School of Business, University of Maryland, College Park, Maryland 20742, USA e-mail: gshmueli@rhsmith.umd.edu.*

Our collaboration with information systems and marketing colleagues has shown just how much the two sides can benefit from crossing the road. eCommerce data are different than other types of data in many, many ways, and they pose real research challenges. Using off-the-shelf statistical methods, which is the current practice, can lead to incorrect or inaccurate results, and, furthermore, to missing out on important real effects. The integration of statistical thinking into the entire process of collecting, cleaning, displaying and analyzing eCommerce data can lead to more sound science and to real research advances.

## 2. FORMAT AND CONTENTS OF THIS SPECIAL ISSUE

This issue has a somewhat different format from standard *Statistical Science* issues for the purpose of capturing the true interdisciplinary nature of eCommerce research. Rather than stand-alone papers, this special issue features mostly "paper-bundles," which combine contributions from authors from the different fields. Each bundle is dedicated to a certain statistical topic, and consists of a methods paper followed by one or more shorter application papers. The methods papers discuss general statistical issues such as data collection and exploration, modeling, methods or inference. The application papers describe a specific eCommerce application within that topic, covering a wide range of state-of-the-art eCommerce research. The goal is to familiarize the statistics community with the types of research questions that are of interest, the types of novel data that arise, methods currently employed, and data-related challenges and questions in eCommerce research.

The topics in this issue touch upon many different statistical areas: Starting with the initial data-collection step, the paper by Bapna et al. describes the process of automated data collection from eCommerce websites and surrounding challenges in data validation, storage and analysis. The bundled application paper by Ghose and Sundararajan describes an effort to quantify pricing strategies for software products using data collected from www.amazon.com.

eCommerce data sets are usually very large, in both dimensionality and size. This requires methods that are suitable for mining this massive amount of information. The data mining bundle includes a paper by Banks and Said, which discusses data mining challenges and opportunities that arise in the eCommerce context. The application paper, by Goldfarb and Lu, describes the fitting of household-level regression models for identifying important factors in Internet portal choice among households.

eCommerce data often arrive as a combination of longitudinal and cross-sectional data, and in many cases they represent an underlying continuous process. The bundle on functional data analysis (FDA) consists of a methodology paper by Jank and Shmueli, that advocates the usefulness of FDA in the eCommerce context for capturing the dynamic nature of the web. The two application papers use FDA to analyze dynamics of prices in online auctions for Indian modern art (Reddy and Dass), and to study the evolution of project complexity in open-source software (Stewart, Darcy and Daniel).

Although designed experiments are possible in the eCommerce context, most data is of an observational nature. Thus, the need for inferring causality is as important in the online setting as it is in the traditional offline setting. The corresponding paper-bundle starts with an overview of propensity scores by Rubin and Waterman. The application paper by Mithas, Almirall and Krishnan illustrates an application of propensity scoring to customer-relationship-management (CRM) systems.

The fifth bundle is devoted to the vibrant area of social networks. The paper by Hill, Provost and Volinsky describes a method for capturing the effect of viral marketing in large networks. The application paper by Dellarocas and Narayan suggests a way to quantify the effect of word-of-mouth on online movie rating using data from Yahoo! Movies.

Finally, we have two unbundled papers, both of which combine methods and applications. The first paper, by Borle, Boatwright and Kadane, takes a Bayesian approach and uses hierarchical models to estimate bidding patterns on eBay. The second paper, by Fienberg, discusses privacy issues associated with the wide availability of electronic data from eCommerce and their possible linkage from multiple sources.

The different papers in this issue give a glimpse into the richness of eCommerce research, and the many opportunities for statisticians to participate in this exciting field. We hope that you will enjoy reading this special issue, and more so, will decide to join this new research field.

## ACKNOWLEDGMENTS

We thank the many referees from statistics and other fields who assisted in greatly improving this special issue. We also thank the enthusiastic workshop participants who have led to the birth of this special issue, and the Executive Editor, Edward George, for his constant support and innovative spirit.

The workshop was supported by National Science Foundation Grant ATS/01-5-24570, the Center of Electronic Markets and Enterprises (`www.smith.umd.edu/ceme/`) and the Statistics Consortium at the University of Maryland (`www.statconsortium.umd.edu`).